\documentclass[12pt]{article} 
\usepackage{amsmath,amsfonts,amscd,a4wide,theorem} 
\theoremheaderfont{\bfseries\rmfamily\upshape} 
                \newtheorem{theorem}{Theorem} 

                \newtheorem{lemma}{Lemma}

{\theorembodyfont{\rmfamily}

                \newtheorem{remark}{Remark}

}

\newcommand{\sech}{\operatorname{sech}} 
\newcommand{\csch}{\operatorname{csch}}

\newcommand{\R}{\mathbb{R}} 
\newcommand{\C}{\mathbb{C}}

\newcommand{\imag}{i}

\newcommand{\e}{{\bf e}_}

 \newcommand{\U}{{\bf U}_}
\newcommand{\w}{{\bf \omega}_}

\newcommand{\qed}{\text{\rule{.4em}{1.7ex}\hspace{.6em}}} 
\newenvironment{proof}[1][]{\noindent {\bf Proof#1:\ }} 
        {\hspace*{.1em}\hfill\qed\bigskip \noindent} 
\newcounter{rom} 
\renewcommand{\therom}{(\roman{rom})} 
{\end{list}} 
\title{Transforms for minimal surfaces in the 5-sphere.}
\begin{document}  
\author{J. Bolton and  L. Vrancken} 

\date{}

\maketitle 

\sloppy 
 
\begin{abstract}\noindent  
We define two transforms between  minimal surfaces with non-circular ellipse of curvature in the 5-sphere, and show how this enables us to construct, from one such surface, a sequence of such surfaces. We also use the transforms to show how to associate to such a surface a corresponding ruled minimal Lagrangian submanifold of complex projective 3-space, which gives the converse of a construction considered in a previous paper, and illustrate this explicitly in the case of bipolar minimal surfaces. 
\end{abstract}

{{\bfseries Key words}:  {\em Sphere, minimal surface, ellipse of curvature,  Lagrangian submanifold, complex projective space}.

{\bfseries Subject class: } 53B25, 53B20.}

\footnotetext{We are grateful to the London Mathematical Society for financial support during the preparation of this paper.}

\section{Introduction} 

Let $f: S \rightarrow S^5(1)$ be a minimal immersion of a surface $S$ into the unit 5-sphere. The image of the unit circle in the tangent space under the second fundamental form of
$f$ is a central planar ellipse $E$ in the normal space  of
$f$ called the {\it ellipse of
curvature}. Let $0\le \theta \le \pi/2$ be such that $\cos \theta$ is the ratio of the lengths
 of the minor and major axes of $E$. The geometrical significance of $\theta$ lies in the fact
 that if  $R_{\theta}$ is the rotation of the  normal space through angle $\theta$ about the
 minor axis of $E$ then  the orthogonal projection  of $R_{\theta}(E)$ onto the plane
 containing $E$ is a circle. If $N$ is the unit vector in the normal space orthogonal to
 the plane containing $E$, then the transform we will consider is obtained by applying 
 $R_{\theta}$ to $N$. Of course, there are certain choices of sign and orientation to be made here,
 and the various choices available give two essentially different transforms.  We shall show that
 these transformed surfaces are also minimal, and that the two transforms are mutual inverses. This enables us to define a sequence $\{f^p \ :\ p\in {\mathbb Z}\}$ of minimal immersions into  $S^5(1)$ with $f^0=f$, and we instigate an investigation of this sequence. 

The transforms  described above are natural generalisations of the polar construction \cite{la} for
 superconformal minimal surfaces in $S^3(1)$ and  $S^5(1)$ (although when the ambient space is  $S^3(1)$ the polar is simply the unit normal to the immersion). In the latter situation, the  ellipse of
curvature is a circle, the angle of rotation $\theta$ is zero, and in both situations the sequence of minimal immersions is periodic with period two. Our motivation for discovering and studying these transforms
 comes from \cite{bovr}, where we showed how to associate two minimal surfaces in $S^5(1)$ to a
 ruled minimal Lagrangian submanifold of complex projective 3-space. We showed that these minimal
 surfaces were related by the above transforms. In the present paper we show that the
 construction described in  \cite{bovr} may be reversed, thus showing that all minimal
 surfaces in $S^5(1)$ whose ellipse of
curvature is  not a circle may be constructed using the methods of  \cite{bovr}. We then illustrate this explicitly in the case of bipolar minimal surfaces \cite{la} in $S^5(1)$. 

\section{Minimal surfaces in  $S^5(1)$}

For the rest of the paper, $f
: S \rightarrow S^5(1)$ 
will denote a minimal immersion
of an oriented surface $S$ into $S^5(1)$.  We use the
orientation and induced metric to give  $S$ the structure of a Riemann surface
in such a way that $f$ is a conformal immersion. If  $II$ denotes the  second fundamental form of
$f$ in $S^5(1)$ we recall that for each $p\in S$ the subset $E(p)$ of the first normal space of $S$ at $p$ given by
$$
E(p)=\{II(X,X) \ |\ X \hbox{ is a unit tangent vector to }S\hbox{ at }p\}
$$
is a (possibly degenerate) central ellipse called the {\it ellipse of
curvature} of $S$ at {p}. 

From now on, we assume that $S$ has non-degenerate non-circular ellipse of 
curvature at every point. In this section we show how to  
associate a complex moving frame to such an immersion, and obtain the moving frame equations and integrability conditions. The approach
we use is based on the theory of harmonic sequences, which is described in \cite{bowo} for the more general situation of minimal surfaces in $S^m(1)$ or $\mathbb CP^m(4)$.

Let $z = x+iy$ be a local complex coordinate
on $S$, and denote $\tfrac{\partial}{\partial z}$ by $\partial$ and
$\tfrac{\partial}{\partial \bar z}$ by $\bar \partial$.
We introduce $\mathbb
C^6$-valued functions $f_0$, $f_1$, $f_2$ by
\begin{align}
&f_0=f,\label{eq0}\\ 
&f_1 = \partial f,\label{eq1}\\
&f_2 = II(\partial, \partial), \label{eq2}
\end{align}
where $II$ now denotes  the complex bilinear extension of the second fundamental form of
$f$ in $S^5(1)$. If $(\quad,\quad)$ is the complex bilinear extension of the standard inner
product on $\mathbb R^6$, it follows that $(f_0,f_1) = 0$ while  conformality of $f$ is equivalent to
\begin{equation}
(f_1,f_1) = 0. \label{eq3}
\end{equation}
Thus  $f_0$, $f_1$, $\bar f_1$ are mutually 
orthogonal and $f_2$ is the component of $\partial f_1$ 
orthogonal to $f_0$, $f_1$, $\bar f_1$.
We note that, by 
Takahashi's Lemma, the minimality of $f$ is equivalent to
$ \partial\bar  \partial f_0=\mu f_0$ for some
$\mu\in \mathbb R$.

If $f_2 = a - ib$ where $a,b$ are $\mathbb R^6$-valued functions, it follows from conformality of $f$  that the ellipse of curvature is homothetic to the image of the map
\begin{equation} \label{ell}
\psi \mapsto II(\cos \psi \tfrac{\partial}{\partial x} + \sin \psi
\tfrac{\partial}{\partial y},
 \cos \psi \tfrac{\partial}{\partial x} + \sin \psi
\tfrac{\partial}{\partial y}) = 2 (a \cos 2\psi +
b \sin 2\psi).
\end{equation}
Since the ellipse of curvature is not a circle we see that 
$(f_2,f_2) \ne 0$. As noticed by Hopf, the function  $(f_2,f_2)$ is holomorphic, so that \cite{bovrwo} there exists a complex coordinate $z$ (which we will call an {\em adapted complex coordinate} for $f$), defined up to rotations by $\tfrac{\pi}{2}$, such that
$(f_2,f_2)=-1$. In this case 
$$-1= (f_2,f_2) =(a,a)-(b,b) -2 i (a,b),$$
so that 
\begin{equation}\label{orth}
(a,a)-(b,b)=-1, \quad (a,b) =0. 
\end{equation}
It now follows from \eqref{ell} that $a$  lies along 
 the minor axis and $b$ the major axis of $E$. It follows from \eqref{orth} that there is a positive function $\phi$ such that 
\begin{equation}
|a|=\sinh\phi, \qquad |b|=\cosh\phi,\label{defphi}
 \end{equation}
 so that the eccentricity $e$ of $E$ is given by
\begin{equation}
 e=\sqrt{1 -\tfrac{\vert a \vert^2} {\vert b \vert^2}} = \sech\phi.
 \end{equation}
We complete our complex moving frame for $f$ by letting $N$ be a real unit vector orthogonal to  $\{f_0,f_1,\bar f_{1},f_2,\bar f_{2}\}$.

It is then straightforward to check that if ${\cal F}=\{f_0,f_1,\bar f_{1},f_2,\bar f_{2},N\}$ and if  $\omega=\log \vert f_1\vert^2$,
then the matrix $A$ of complex bilinear inner products of the frame vectors of $\cal F$ is
given by
\begin{equation}\label{matA}
A=\begin{pmatrix}
&1&0&0&0&0&0\\
&0&0&e^\omega&0&0&0\\
&0&e^\omega&0&0&0&0\\
&0&0&0&-1&\cosh 2\phi&0\\
&0&0&0&\cosh 2\phi&-1&0\\
&0&0&0&0&0&1\end{pmatrix}.\end{equation}

We now write down the moving frame equations for $\cal F$.
A straightforward computation using \eqref{matA} shows that if  $\alpha = (\partial f_2,N)$,  
then the moving frame equations for $\cal F$ may be written in terms of $\omega$, $\phi$, $\alpha$ as follows:
\begin{align*}  
&\partial f_0=f_1,\\
&\partial f_1 =f_2 + \partial \omega\, f_1,\\
&\partial \bar f_1=- e^{\omega} f_0,\\
&\partial f_2=e^{-\omega}\bar f_{1} +2 \partial \phi\,\coth2\phi\ f_2 + 2 \partial \phi\,\csch2\phi\  \bar f_{2} +\alpha N,\\
&\partial \bar f_{2}=-e^{-\omega}\cosh2\phi \ \bar f_1,\\
&\partial N=-\alpha \csch^22\phi\ (f_2 +\cosh2\phi \, \bar f_{2}).
\end{align*}
Of course, the corresponding $\bar \partial$ equations may be found by taking the conjugates of the above. It follows from uniqueness of solutions of linear differential equations and the integrability 
conditions $\partial\bar\partial{\cal F}=  \bar\partial\partial{\cal F}$ that  a minimal surface with non-circular non-degenerate ellipse
of curvature in $S^5(1)$ is determined, up to $O(6)$-congruence, 
by functions $\omega$, $\phi$, $\alpha$ satisfying the following system of differential equations:
 \begin{equation}\label{intF}\begin{aligned}
 &\bar \partial \alpha = -2  \bar\alpha \partial \phi\,\csch2\phi,\\
 &\bar \partial \partial \omega = -e^\omega+e^{-\omega} \cosh2\phi,\\
 &2\bar \partial \partial \phi = \alpha \bar \alpha\csch2\phi- 
e^{-\omega}\sinh2\phi. 
 \end{aligned}\end{equation} 
The functions $\alpha$, $\omega$ and $\phi$ have geometrical significance:  $\alpha$ is analogous to the torsion of a space curve in that it is a measure of the rate at which the surface is pulling away from the great 4-sphere which contains its tangent and first normal space, the metric on the surface is given by $2e^\omega\vert dz\vert^2$, and  $\phi$ is a measure of the eccentricity of the ellipse of curvature.
\begin{remark}
The above equation for $\phi$ may be used to show that every compact minimal surface $S$ in $S^4(1)$ contains at least one point at which the ellipse of curvature is a point, a line segment or a circle. Otherwise,  $\phi$ would be  a smooth globally defined positive function on $S$ satisfying the two-dimensional sinh-Gordon equation.
\end{remark}

\section{The transforms}

In this section we show how to associate to $f$ two other minimal immersions of the Riemann surface $S$ into $S^5(1)$ which both induce the same conformal strucure on $S$ as that induced by $f$, and also have non-degenerate non-circular ellipse of 
curvature at every point. We further show that an adapted complex coordinate $z$ for $S$ is also an adapted complex coordinate for our new minimal immersions. 
 Recall from the introduction that these minimal immersions are obtained by rotating $N$ in the normal space through a geometrically significant angle $\theta$ about the minor axis of the ellipse of curvature. So, let $z$ be an adapted complex coordinate for $f$ and let $\cos\theta=\vert a \vert /\vert b \vert= \tanh\phi$. Then, applying rotations of $\pm \theta$ to $\pm N$ gives the  four possibilities
$$
\pm\tfrac{1}{\cosh \phi} \tfrac{b}{\vert b \vert } \pm\tanh \phi N.
$$

For definiteness, we take $N$ to be such that $\{f_0,\tfrac{\partial f_0}{\partial x},\tfrac{\partial f_0}{\partial
y},II( \tfrac{\partial f_0}{\partial x},\tfrac{\partial f_0}{\partial x}),II(\tfrac{\partial f_0}{\partial
x},\tfrac{\partial f_0}{\partial y}),N\}$ is  a positively oriented orthogonal moving frame of ${\mathbb R}^6$, and  the two transforms we will consider are those
 given by 
\begin{align}
f^+&=-\tfrac{1}{\cosh \phi} \tfrac{b}{\vert b \vert } +\tanh \phi N,\label{plus}\\
f^-&=-\tfrac{1}{\cosh \phi} \tfrac{b}{\vert b \vert } -\tanh \phi N.\label{min}
\end{align}

Thus if we orient the normal space by taking $\{II(\tfrac{\partial f_0}{\partial
x},\tfrac{\partial f_0}{\partial x}),II(\tfrac{\partial f_0}{\partial
x},\tfrac{\partial f_0}{\partial y}),N\}$ to be positively oriented, then the {\em $(+)$transform} $f^+$ is obtained from $f$ by the {\em $(+)$construction} which is given by \eqref{plus} and consists of rotating $N$ about the minor axis of the ellipse of curvature through the angle $\theta$ anticlockwise, while  the {\em $(-)$transform} $f^-$ is obtained  from $f$ by the {\em $(-)$construction}  which  is given by \eqref{min} and consists of rotating $-N$ about the minor axis through the angle $\theta$ clockwise. 

We note the following for later use. Let vol denote the complexification of the standard volume form of ${\mathbb R}^6$. Since $\det A= e^{2\omega}\sinh^22\phi$, we see that vol$(f_0,f_1,\bar f_{1},f_2,\bar f_{2},N)=\pm  e^{\omega}\sinh 2\phi$.
However,  
$$
\hbox{vol}(f_0,f_1,\bar f_{1},f_2,\bar f_{2},N)=-\tfrac{1}{4} \hbox{vol}\left(f_0,\tfrac{\partial f_0}{\partial x},\tfrac{\partial f_0}{\partial
y},II( \tfrac{\partial f_0}{\partial x},\tfrac{\partial f_0}{\partial x}),II(\tfrac{\partial f_0}{\partial
x},\tfrac{\partial f_0}{\partial y}),N\right),
$$ 
so that
 \begin{equation}\label{volume1}
\hbox{vol}(f_0,f_1,\bar f_{1},f_2,\bar f_{2},N)=-  e^{\omega}\sinh 2\phi.
\end{equation}

We now show that $f^+$ and $f^-$ both induce the same conformal structure on $S$ as that induced by $f$. In order to do this, we first define
\begin{equation*}
f_1^\epsilon =\partial f^\epsilon.
\end{equation*}
Here and subsequently, we use $\epsilon$ as a superscript taking value $+$ or $-$, and use $\epsilon =\pm 1$ in the corresponding equations as appropriate.
Then, using \eqref{plus}, \eqref{min} and the moving frame equations for $\cal F$, we find that 
\begin{equation}
f_1^\epsilon =-i e^{-\omega}\bar f_{1} - \tfrac{1}{2}(\alpha \epsilon +2 i \partial \phi) \sech^2 \! \phi \ \Big(\csch2\phi\, f_2 + \coth 2\phi\, \bar f_{2} +\epsilon i N\Big).\label{dfeps}
\end{equation}

From this, a computation using \eqref{matA} shows that 
\begin{align}
&(f_1^\epsilon,f_1^\epsilon)=0,\label{isofe}\\
&\vert f_1^\epsilon  \vert^2=e^{-\omega}+ \tfrac{1}{2} \vert  \alpha \epsilon +2 i \partial
\phi \vert^2\sech^2\!\phi.\label{omeps} 
\end{align}
Thus the maps $f^\epsilon$ define conformal immersions of $S$ into
$S^5(1)$, so that if $z=x+iy$ then $(x,y)$ are isothermal coordinates not only for the original immersion $f$ 
but also for the two newly constructed immersions $f^\epsilon$. We note that the metric induced on $S$ by $f^\epsilon$ is given by $2e^{\omega^\epsilon}|dz^2|$, where $\omega^\epsilon=\log\vert f_1^\epsilon  \vert^2$. 

We now show that each $f^\epsilon$ is minimal. By Takahashi's lemma, it is sufficient to check that
$\partial \bar \partial f^\epsilon$ is a multiple of $f^\epsilon$. We first note that it follows quickly from \eqref{dfeps} that 
\begin{align*}
&(\partial \bar \partial f^\epsilon,f_0)=0,\\
&(\partial \bar \partial f^\epsilon,f_1)=0,
\end{align*}
while a straightforward 
computer calculation shows that
\begin{align*}
&(\partial \bar \partial f^\epsilon,f_2)=-i(e^{-\omega}+ \tfrac{1}{2}  \vert  \alpha \epsilon +2 i \partial
\phi \vert^2\sech^2\!\phi), \\
&(\partial \bar \partial f^\epsilon,N)=-\epsilon  \tanh\phi\ (e^{-\omega}+ \tfrac{1}{2}  \vert  \alpha \epsilon +2 i \partial
\phi \vert^2\sech^2\!\phi).
\end{align*}
However, it follows from \eqref{defphi}, \eqref{matA}, \eqref{plus} and \eqref{min} that
\begin{align}
&(f^\epsilon,f_0)=0,\label{ortho1}\\
&(f^\epsilon,f_1)=0,\label{ortho2}\\
&(f^\epsilon,f_2)=i,\label{ortho3}\\
&(f^\epsilon,N)=\epsilon \tanh\phi\label{ortho4},
\end{align}
and hence, using \eqref{omeps},
\begin{equation*}
\partial \bar \partial f^\epsilon=-\vert f_1^\epsilon  \vert^2 f^\epsilon,
\end{equation*}
showing that $\partial \bar \partial f^\epsilon$ is indeed a multiple of $f^\epsilon$, so that each $f^\epsilon$ is a minimal immersion of $S$ into 
$S^5(1)$. 

Finally in this section, we show that $z$ is also an adapted complex coordinate for $f^\epsilon$. If $II^\epsilon$ denotes the second fundamental form of $f^\epsilon$, we put
$$
 f_2^\epsilon=II^\epsilon(\partial,\partial)=a^\epsilon -i b^\epsilon.
$$

As in the moving frame equations for $\cal F$, if $\omega^\epsilon=\log\vert f_1^\epsilon \vert^2$, then
\begin{equation}
f_2^\epsilon = \partial f_1^\epsilon -\partial \omega^\epsilon f_1^\epsilon.\label{f2eps}
\end{equation}
It follows from \eqref{plus}, \eqref{min} and the moving frame equations for $\cal F$ that  
\begin{equation}\label{df1eps}
\begin{split}
\partial f_1^\epsilon &= i f_0 + e^{-\omega}\Big(i \partial \omega +\tanh \phi\  (\alpha \epsilon + 2i
\partial \phi)\Big) \bar f_1+ 2 \nu \csch^4\!2\phi\ \sinh^2\!\phi\  f_2 \\
&\qquad + (1/2)\nu \coth2 \phi \csch2\phi \sech^2\!\phi \ \bar f_{2} 
  +i \epsilon \nu\csch^22\phi \tanh \phi\ N,
\end{split}
\end{equation} 
where 
\begin{equation*}
\begin{split}
\nu &=2 \alpha \epsilon\partial \phi (-2 +\cosh2\phi ) +8 i \sinh^2\!\phi\ (\partial \phi)^2\\
&\qquad - \epsilon \partial \alpha\ \sinh 2\phi  +i \alpha^2 -2 i \sinh 2\phi\ 
\partial \partial \phi,
\end{split}
\end{equation*}
so that, using \eqref{matA}, 
\begin{equation*}
(\partial f_1^\epsilon,\partial f_1^\epsilon)=-1.
\end{equation*}
Equations \eqref{isofe} and \eqref{f2eps} now imply that
$(f_2^\epsilon,f_2^\epsilon)=-1$, so that $z$ is also
an adapted complex coordinate for each $f^\epsilon$.  In particular, each $f^\epsilon$ has  non-degenerate non-circular ellipse of curvature.

Summarizing the above, we have the following theorem:
\begin{theorem} 
Let $f:S \rightarrow S^5(1)$ be a minimal immersion with non-degenerate non-circular 
ellipse of curvature at every point. Then the $(+)$transform $f^+$  and the $(-)$transform $f^-$ of $f$ are both
minimal immersions of $S$ into $S^5(1)$ which induce 
 the same conformal structure on $S$ as that induced by $f$. Moreover,  both  $f^+$ and $f^-$ have  non-degenerate non-circular ellipse of curvature at every point and an adapted complex coordinate for $f$ is also an  adapted complex coordinate for $f^+$ and $f^-$.
\end{theorem}

\begin{remark} The final statement of the above theorem is equivalent to saying that $f^+$ and $f^-$ both have the same the 
have the same $U_{2,-2}$ invariant (see \cite{bowo} for the definition of this and related invariants).
\end{remark}

\section{A symmetric adapted moving frame}
In this section we begin the study of $f$ and its transforms by constructing a moving frame $\mathcal B$ which gives equal prominence to $f$ and $f^{\epsilon}$. We also obtain the moving frame equations and integrability conditions for $\mathcal B$. 

So, let ${\mathcal B}=\{f_0,f_1,\bar{f_1},\bar{f_1^\epsilon}, f_1^\epsilon, f_0^\epsilon\}$. 
It follows quickly from \eqref{matA} and the moving frame equations for $\mathcal F$ that the
matrix $B$ of inner products of elements of  ${\mathcal B}$ is given by
\begin{equation}\label{matb}
B = \begin{pmatrix}
&1&0&0&0&0&0\\
&0&0&e^\omega&0&-i&0\\
&0&e^\omega&0&i&0&0\\
&0&0&i&0&e^{\omega^\epsilon}&0\\
&0&-i&0&e^{\omega^\epsilon}&0&0\\
&0&0&0&0&0&1
\end{pmatrix},
\end{equation} 
from which we see that 
\begin{equation}\label{detB}\det B= (e^{\omega+\omega^\epsilon}-1)^2,
\end{equation}
 implying that the vectors in  ${\cal B}$ are linearly independent as long as $\omega+\omega^\epsilon\ne 0$. 

\begin{lemma} \label{lem1} We have that $\omega + \omega^\epsilon > 0$ on an open dense
subset of $S$. 
\end{lemma}
\begin{proof} 
It follows from \eqref{omeps} that $\omega +
\omega^\epsilon \ge  0$, and that $\omega + \omega^\epsilon = 0$ 
on an open set 
if and only if
$\alpha=-2 i \epsilon \partial
\phi$. 
Taking the derivative of this expression with respect to $\bar \partial$
and using the integrability conditions \eqref{intF} for the moving frame equations for $\cal F$ it then
follows that
\begin{equation*}
-4i  \epsilon \csch2 \phi\ \partial \phi  \bar \partial \phi = 
e^{-\omega}  i\epsilon \csch2\phi\ (\sinh^2\!2 \phi - 4 e^\omega \partial \phi \bar \partial\phi).
\end{equation*}
Simplifying the above equation then yields the contradiction $0=e^{-\omega}\sinh 2\phi$.
\end{proof}

The advantage of the above condition is that on this open dense subset we can 
investigate the original immersion $f$ and the new immersion $f^\epsilon$  with 
respect to the frame ${\cal B}=\{f_0,f_1,\bar{f_1},\bar{f_1^\epsilon}, f_1^\epsilon, f_0^\epsilon\}$.

It follows from  \eqref{detB} that vol$(f_0,f_1,\bar{f_1},\bar{f_1^\epsilon}, f_1^\epsilon, f_0^\epsilon)= \pm (e^{\omega+\omega^\epsilon}-1)$. In order to determine the sign we compute the volume explicitly. 
A straightforward calculation using \eqref{volume1} and \eqref{dfeps} yields that
\begin{align*}
\hbox{vol}(f_0,f_1,\bar{f_1},\bar{f_1^\epsilon}, f_1^\epsilon, f^\epsilon)&=\epsilon \csch^2\!2\phi\ \tanh\phi \ \vert \alpha \epsilon +2 i \partial
\phi \vert ^2  \hbox{vol }(f_0,f_1,\bar f_{1},f_2,\bar f_{2},N) \\
&=-\epsilon e^\omega \csch 2\phi\ \tanh \phi \ \vert \alpha \epsilon +2 i \partial
\phi \vert ^2\\
&=-\epsilon \tfrac{1}{2} \sech^2 \!\phi\  e^\omega \vert  \alpha \epsilon +2 i \partial
\phi \vert^2.
\end{align*}
Thus, using \eqref{omeps},
 \begin{equation}\label{volume2}
\hbox{vol}(f_0,f_1,\bar{f_1},\bar{f_1^\epsilon}, f_1^\epsilon, f_0^\epsilon)
=-\epsilon (e^{\omega+\omega^\epsilon}-1).
\end{equation}

We now intoduce the function $\gamma^\epsilon$ by
\begin{equation}\label{gaminv}
\gamma^\epsilon=(\partial f_1,f_1^\epsilon),
\end{equation} 
and  use  \eqref{matb} to write down the moving frame equations
 for $\cal B$ in terms of $\omega$, $\omega^\epsilon$, $\gamma^\epsilon$ as follows. 
\begin{align*}  
&\partial f_0=f_1,\\
&\partial f_1 =\tfrac{-i\gamma^\epsilon+\partial\omega e^{\omega+\omega^\epsilon}}
{e^{\omega+\omega^\epsilon}-1}f_1 + \tfrac{e^{\omega}(\gamma^\epsilon +
i\partial\omega)}{e^{\omega+\omega^\epsilon}-1}\bar f_1^\epsilon + if_0^\epsilon,\\
&\partial \bar f_1= -e^{\omega}f_0,\\
&\partial \bar f_{1}^\epsilon=-e^{\omega^\epsilon} f_0^\epsilon,\\
&\partial f_1^\epsilon=if_0 - \tfrac{e^{\omega^\epsilon}(\gamma^\epsilon -
i\partial\omega^\epsilon)}{e^{\omega+\omega^\epsilon}-1}\bar f_1 + 
\tfrac{i\gamma^\epsilon+\partial\omega^\epsilon e^{\omega+\omega^\epsilon}}
{e^{\omega+\omega^\epsilon}-1}f_1^\epsilon,\\
&\partial f_0^\epsilon=f_1^\epsilon.
\end{align*}
As before, the corresponding $\bar \partial$ equations may be found by taking the conjugates of the above. The integrability 
conditions of the above system of equations are
 \begin{equation}\label{eq24} \begin{aligned}
 &\bar \partial \gamma^\epsilon=i(e^{\omega}-e^{\omega^\epsilon}),\\
 &\partial\bar{\partial}\omega=-2 \sinh \omega
 +\tfrac{1}{e^{\omega+\omega^\epsilon}-1} \vert \gamma^\epsilon+i \partial
 \omega\vert^2,\\
 &\partial\bar{\partial}\omega^\epsilon= -2 \sinh \omega^\epsilon 
 +\tfrac{1}{e^{\omega+\omega^\epsilon}-1} \vert \gamma^\epsilon-i \partial
 \omega^\epsilon\vert^2.
\end{aligned}\end{equation}
As before, solutions of \eqref{eq24} correspond up to $O(6)$ congruence to a minimal surface and its 
$\epsilon$-transform. 

\section{The (+) and ($-$) constructions are mutual inverses}

In this section we use the frame ${\cal B}$ introduced in the previous section to prove the following.

\begin{theorem} Let $f:S \rightarrow S^5(1)$ be a minimal immersion with non-degenerate non-circular 
ellipse of curvature at every point. Then the  
 $(+)$construction  and the $(-)$construction are mutual  inverses in the sense that both  $f^+$ and $f^-$ have non-degenerate non-circular ellipse of curvature and  $(f^+)^-=(f^-)^+=f$.
\end{theorem}
\begin{proof} We denote by
$p_{\tilde\epsilon}(f_0^\epsilon)$, where $\tilde \epsilon \epsilon=-1$, the
image of $f_0^\epsilon$ by the $\tilde \epsilon$ construction.
So, by \eqref{plus} and \eqref{min}, we have that
\begin{equation}
p_{\tilde\epsilon}(f_0^\epsilon)=-\tfrac{1}{\cosh \phi^\epsilon} 
 \tfrac{b^\epsilon}{\vert b^\epsilon \vert } +\tilde \epsilon 
\tanh \phi^\epsilon  N^\epsilon = -\tfrac{b^\epsilon}{\cosh^2 \!\phi^\epsilon} 
 +\tilde \epsilon 
\tanh \phi^\epsilon  N^\epsilon.
\end{equation}

Now let 
\begin{equation*}
v=\tfrac{b^\epsilon}{\cosh^2\! \phi^\epsilon} +f_0.
\end{equation*}
It follows from the $\epsilon$-analogue of \eqref{matA} and from \eqref{matb} that $v$ is orthogonal to $f_0^\epsilon$ and $f_1^\epsilon$, while, using the  $\epsilon$-analogues of \eqref{orth} and \eqref{defphi}, we see that
\begin{align*}
(v,f_2^\epsilon)& = (\tfrac{b^\epsilon}{\cosh^2 \phi^\epsilon},f_2^\epsilon)+(f_0,f_2^\epsilon)\\
&=-i\tfrac{\vert b^\epsilon\vert^2}{\cosh^2 \phi^\epsilon}+(f_0,\partial f_1^\epsilon)\\
&=-i-(\partial f_0,f_1^\epsilon)\\
&=-i-(f_1,  f_1^\epsilon)\\
&=-i+i=0,
\end{align*}
so that $v$ is a scalar multiple of $N^\epsilon$. We now use the 
moving frame equations for $\cal B$ to determine this scalar multiple as follows. First note that 
\begin{align*} 
 \hbox{vol}(f_0^\epsilon,f_1^\epsilon,
 \bar f_{1}^\epsilon,f_2^\epsilon, \bar f_{2}^\epsilon,v) 
&= \hbox{vol}(f_0^\epsilon,f_1^\epsilon,
 \bar f_{1}^\epsilon,f_2^\epsilon, \bar f_{2}^\epsilon,f_0)\\
&=\hbox{vol}(f_0^\epsilon,f_1^\epsilon,
 \bar f_{1}^\epsilon,\partial f_1^\epsilon, \bar \partial \bar
 f_{1}^\epsilon,f_0)\\
&= -\tfrac{e^{2\omega^\epsilon}}{(e^{\omega+\omega^\epsilon}-1)^2}  
\vert \gamma^\epsilon - i \partial\omega^\epsilon \vert^2  \hbox{vol}(f_0,f_1,\bar{f_1},\bar{f_1^\epsilon}, f_1^\epsilon, f^\epsilon)\\
&=\epsilon
\tfrac{e^{2\omega^{\epsilon}}}{(e^{\omega+\omega^\epsilon}-1)}\vert \gamma^\epsilon 
- i \partial\omega^\epsilon \vert^2,
 \end{align*}
where the final  equality above is obtained using \eqref{volume2}.

However, it follows from  \eqref{eq24} and the 
$\epsilon$-analogue of \eqref{intF} that 
\begin{equation}\label{sinhphi}
\tfrac{e^{\omega^{\epsilon}}}{(e^{\omega+\omega^\epsilon}-1)}\vert \gamma^\epsilon- i \partial\omega^\epsilon \vert^2=2 \sinh^2 \!\phi^\epsilon, 
\end{equation}
so that 
\begin{equation*}
 \hbox{vol}(f_0^\epsilon,f_1^\epsilon,
 \bar f_{1}^\epsilon,f_2^\epsilon, \bar f_{2}^\epsilon,v) = 2\epsilon e^{\omega^\epsilon}\sinh^2 \!\phi^\epsilon.
\end{equation*}

We next note that the $\epsilon$-analogue of \eqref{volume1} gives
 \begin{equation*}
 \hbox{vol}(f_0^\epsilon,f_1^\epsilon,
 \bar f_{1}^\epsilon,f_2^\epsilon, \bar f_{2}^\epsilon,N^\epsilon) =-e^{\omega^\epsilon}\sinh 2\phi^\epsilon,
 \end{equation*}
 so it follows that 
 $
 v=-\epsilon \tanh \phi^\epsilon N^\epsilon. 
$

Thus, from (24) we have that 
\begin{align*}
p_{\tilde\epsilon}(f_0^\epsilon)&=-v+f_0+\tilde \epsilon 
\tanh \phi^\epsilon  N^\epsilon \\
&=(\tilde \epsilon +  \epsilon)\tanh \phi^\epsilon  N^\epsilon +f_0=f_0,
\end{align*} 
implying that the $(+)$construction and the $(-)$construction are mutual inverses.
\end{proof}

The above theorem shows that we may associate to a minimal  immersion $f:S\to S^5(1)$ with 
non-degenerate non-circular ellipse of curvature a sequence $\{f^p : p\in {\mathbb Z}\}$ of such minimal immersions with $f^0=f$ and, for each $p$,  $f^{p+1}=(f^p)^+$ and $f^{p-1}=(f^p)^-$.  Moreover, each element of the sequence induces the same conformal structure on $S$, and an adapted complex coordinate $z$ for any $f^p$ is an adapted complex coordinate for each element of the sequence. 

\section{The geometry of the invariants}

In the previous section we showed that we may associate to a minimal immersion $f:S \rightarrow S^5(1)$ with non-degenerate
non-circular ellipse of curvature a sequence $\{f^p : p\in {\mathbb Z}\}$ of such minimal immersions with $f^0=f$. For the remainder of the paper we will use the superfix $^p$ to denote objects connected with $f^p$. For instance, with each 
immersion $f^p$ we associate as before the invariants $\omega^p$, $\phi^p$ and $\alpha^p$. Moreover 
with each $(+)$construction, $f^p \rightarrow f^{p+1}$, we associate the invariants $\omega^p$, 
$\omega^{p+1}$ and $\gamma^{p+1}=(\partial f^p_1, f_1^{p+1})$, while with each $(-)$construction
we associate the invariants $\omega^p$, 
$\omega^{p-1}$ and $\delta^{p-1}=(\partial f^p_1, f_1^{p-1})$. Thus, $\gamma^1$
is equal to the invariant $\gamma^+$ used in previous sections,  while $\delta^{-1}$ is equal to $\gamma^-$. Since, from \eqref{matb}, we have that $(f^p_1, f_1^{p+1})=-i$, it is clear that $\delta^{p} = -\gamma^{p+1}$.

As already mentioned, the geometrical meaning of the invariants $\omega^{p}$ is clear, since they give the metric induced on $S$ by $f^p$. Also, 
the final moving frame equation for ${\mathcal F}$ implies that $\alpha^p=0$ if and only if  $f^p$ is not linearly full. In this section we look more closely at this situation, and also obtain a geometrical characterisation of the condition $\gamma^{p+1}=0$.

We begin with a useful lemma.
\begin{lemma}\label{lem2} Let $A$ be an orientation reversing isometry of ${\mathbb R}^6$. Then 
$$
(Af)^-=A(f^+)\ ,\qquad (Af)^+=A(f^-)\ .
$$
In fact, more generally, for each integer $p$ we have that $(Af)^p=A(f^{-p})$.
\end{lemma}
\begin{proof} As $A$ is an orientation reversing isometry, if $N$ is the normal vector associated to $f$ as in equations \eqref{plus} and \eqref{min}, then the corresponding normal associated to $Af$ is $-AN$.  The first result is now clear from the definitions of the $(+)$ and $(-)$constructions given in \eqref{plus} and \eqref{min}. The second may be proved in a similar manner, and the final statement follows by induction.
\end{proof}

\begin{theorem}\label{th3} Let $\{f^p\}$ be the sequence of minimal immersions into $S^5(1)$ determined by a minimal immersion $f=f^0$ with non-degenerate
non-circular ellipse of curvature in $S^5(1)$. For each integer $q$, the following three statements are equivalent:
\begin{enumerate} 
\item  $\alpha^q=0$,
\item $f^q$ is not linearly full,
\item there exists an orientation reversing isometry $A \in O(6)$ such that $f^{q+1} =A(f^{q-1})$.
\end{enumerate}
Moreover, in this case,  for every integer $r$ we have that $f^{q-r}$ and $f^{q+r}$ are  congruent via reflection in the great 4-sphere containing $f^q$.
\end{theorem}
\begin{proof}
We have already noted the equivalence of the first two statements. Now suppose that condition 3 holds. Taking the $(-)$construction of this, we may use Lemma \ref{lem2} to see that $f^q=A(f^q)$. Since $A$ has at least one eigenvalue equal to $-1$, it now follows that $f^q$ is not linearly full and that  $A$ is reflection in the great 4-sphere containing $f^q$.

Conversely suppose that $f^q$ is contained is a totally geodesic $S^4(1)$. In this case $N^q$ is a constant vector, so it is clear from \eqref{plus} and \eqref{min} that $f^{q-1}$ and $f^{q+1}$ are congruent via reflection in the totally geodesic $S^4(1)$ containing $f^q$.  The final statement of the theorem now follows from Lemma \ref{lem2}.
\end{proof}

A similar characterisation also exists for $\gamma^{q+1}=0$. We show in Section 9 that this situation can actually arise; in this case $f^q$ is a bipolar surface in the sense of Lawson \cite{la}.

\begin{theorem} \label{th4}  Let $\{f^p\}$ be the sequence of minimal immersions into $S^5(1)$ determined by a minimal immersion $f=f^0$ with non-degenerate
non-circular ellipse of curvature in $S^5(1)$. For each integer $q$, $\gamma^{q+1}=0$ if and only if there exists an orientation reversing isometry $A \in O(6)$ such that $f^{q+1} =A(f^{q})$.
Moreover, in this case, $A$ is reflection in a great subsphere of $S^6(1)$ and  for every integer $r$ we have that $f^{q+1+r}= A(f^{q-r})$.
\end{theorem}
\begin{proof}
Assume that  $f^{q+1} =A(f^{q})$. Then,
\begin{align*}
\gamma^{q+1}&= (\partial f_1^q,f_1^{q+1})\\
&=(\partial A(f_1^q), A(f_1^{q+1})).
\end{align*}

However, 
\begin{equation*} A(f_1^q)= A\partial f^q = \partial A(f^q)=\partial f^{q+1}=f^{q+1}_1,
\end{equation*}
while, using Lemma \ref{lem2},
\begin{equation*}  A(f_1^{q+1})=A(\partial f^{q+1}) = \partial  A(f^{q+1})=\partial ((A(f^q))^-)= \partial ((f^{q+1})^-)=\partial f^q=f_1^q.
\end{equation*}

Thus 
\begin{equation*}
\gamma^{q+1}= (\partial f^{q+1}_1,f_1^q)= -(f^{q+1}_1 , \partial f_1^q)= -\gamma^{q+1},
\end{equation*}
so that  $\gamma^{q+1}=0$.

Conversely suppose that $\gamma^{q+1} = 0$. It then follows from the integrability conditions \eqref{eq24} that $\omega^q = \omega^{q+1}$. Since the set $\mathcal B$ is a basis for ${\mathbb C}^6$ we may define, for each $z$,  a unique linear map $A(z)$ by
\begin{align*}
&A f_0^q = f_0^{q+1} ,\qquad \qquad\qquad A f_1^q = f_1^{q+1}, \qquad \qquad\qquad A \bar f_1^q = \bar f_1^{q+1},\\
&A f_0^{q+1} = f_0^q ,\qquad \qquad\qquad A f_1^{q+1} = f_1^q, \qquad \qquad\qquad A \bar f_1^{q+1} = \bar f_1^q.
\end{align*}
However, it follows from the moving frame equations for $\mathcal B$ that  $A$ does not depend on $z$, while \eqref{matb} may be used to show that  $A$ is an isometry. It is clear from the definition that $A$ has determinant $-1$ and $A^2$ is the identity, so it follows that $A$ is a reflection. As in the previous theorem, the final statement follows from Lemma \ref{lem2}.
\end{proof}

In the previous two theorems we have considered two situations which led to the conclusion that two elements in a sequence $\{f^p\}$ are congruent via an orientation reversing isometry. An easy application of Lemma \ref{lem2} quickly leads to the following theorem, which shows that the above situations are the only ones for which this can happen.

\begin{theorem} Let $\{f^p\}$ be the sequence of minimal immersions into $S^5(1)$ determined by a minimal immersion $f=f^0$ with non-degenerate
non-circular ellipse of curvature in $S^5(1)$. Suppose that two elements $f^q$ and $f^r$ of the sequence are congruent via an orientation reversing isometry. Then there are two possibilities, depending on the parity of $q-r$. Either
\begin{enumerate} 
\item there exists an integer $s$ for which $\alpha^s=0$, or
\item there exists an integer $s$ for which $\gamma^{s+1} =0$.
\end{enumerate}
\end{theorem}

Finally in this section, we note that if two elements of a sequence $\{f^p\}$ are related by an orientation preserving isometry then the sequence is periodic in a natural sense.  We do not know if this situation can actually arise.  

\section{Minimal surfaces and ruled Lagrangian submanifolds}

In a previous paper \cite{bovr}, we studied minimal  Lagrangian  submanifolds of complex projective 
3-space $\C P^3(4)$ which admit a foliation by asymptotic curves. Such submanifolds
can be divided in to three types.
\begin{enumerate}
\item Those which additionally satisfy Chen's equality \cite {ch}. These were studied and characterized in  \cite{boscvrwo}, \cite{boscvr}, and are closely related to minimal immersions of surfaces in $S^5(1)$ with ellipse of curvature a circle.
\item Those for which the unit tangent vectors to the asymptotic curves form a Killing vector field. It is shown in \cite{bovr} that these are related to minimal surfaces in $S^3(1)$, and classification theorems are obtained in \cite{cavr}.
\item All the rest. In \cite{bovr} we showed how to construct, starting from such a submanifold a pair of minimal immersions of a surface $S$ into $S^5(1)$ with non-degenerate non-circular ellipse of curvature  which are related by the $(+)$ and $(-)$constructions. 
\end{enumerate}
In this section we deal with the converse of the construction described in \cite{bovr}. We show how to associate to a minimal surface in $S^5(1)$ with  non-degenerate non-circular ellipse of curvature a Lagrangian submanifold of $CP^3(4)$ belonging to the last type. We will use the notation and terminology of  \cite{bovr}.

In order to simplify the notation, we denote the immersion of $S$ into $S^5(1)$ by $f$ and let $g$ be the $(+)$transform of  $f$. As usual, we let $z$ be an adapted complex coordinate,  put $f_0=f$, $f_1=\partial f$, $g_0=g$ and $g_1=\partial g$. 
We will show that for a suitable interval $I$ of real numbers,  $M=I\times S$ may be realised as a  minimal  Lagrangian  submanifold of type 3 such that if $\{\e1,\e2,\e3\}$ is an orthonormal frame along $M$ of the type described in \cite{bovr} with $\e1$ being tangential to the asymptotic curves and $\e2,\e3$ being eigenvectors of the second fundamental form $A_{J\e1}$ of $M$ with respect to the normal $J\e1$ with corresponding eigenvalues $\pm \lambda$ ($\lambda>0$), then the corresponding map ${\cal U}=(\U1, \ldots ,\U6):M\to SO(6)$ has
\begin{equation} 
\U2(t,z)=g_0(z), \quad \U4(t,z)=f_0(z),\label{U24}
\end{equation}
where $t$ is the standard real coordinate on $I$.

In fact, we use the invariants $\omega$, $\omega^+$ and $\gamma^+$ to construct a map ${\cal U}=(\U1, \ldots ,\U6):M\to SO(6)$ satisfying \eqref{U24}, with the property that if $\Omega={\cal U}^{-1}d{\cal U}$ then $\Omega$ has the form of (33) of  \cite{bovr} for suitable functions $z_{21}^2$, $z_{12}^3$, $z_{22}^3$, $z_{32}^3$, $\lambda$, $a$ and $b$ on $M$, and linearly independent  1-forms $\omega_1, \omega_2, \omega_3$ on $M$. Having done this, it is straightforward to deduce that we may reverse the construction given in \cite{bovr} in order to construct from such a map $\cal U$ our required Lagrangian submanifold of $CP^3(4)$, with the orthonormal basis  $\{\e1$, $\e2$, $\e3\}$ being the basis of vectors dual to $\{\omega_1, \omega_2, \omega_3\}$.

We begin by noting from (42) of \cite{bovr} that in order for \eqref{U24} to hold we require that
\begin{equation} 
d\U2(\e2-i\e3)=2\theta_1\sqrt \lambda \, g_1,\label{XiY}
\end{equation}
where $\theta_1$ is a fourth root of unity. In fact, we may assume that $\theta_1=1$ by rotating our adapted complex coordinate $z$ through a suitable multiple of $\pi/2$. In a similar manner, 
\begin{equation} 
d\U4(\e2-i\e3)=2\theta_2\sqrt \lambda \, f_1,\label{tXitY}
\end{equation}
for some fourth root of unity $\theta_2$.

We next note that for  (34) and (35) of  \cite{bovr} to hold we must have that 

$$d\U2(\e2-i \e3)= (z_{12}^3-1 -i z_{21}^2)(\U1 +i \U3) - \lambda (\U5+i \U6),$$
and 
$$d\U4({\e2 -i \e3}) =  i\lambda (\U1-i \U3) + i (1+z_{12}^3 -iz_{21}^2) (\U5-i\U6),$$
so that, from \eqref{XiY} and \eqref{tXitY}, 
\begin{align}
2\sqrt{\lambda}\, g_1 &= (z_{12}^3-1 -i z_{21}^2)(\U1 +i \U3) - 
\lambda (\U5+i \U6),\label{reveq1}\\
2   \theta_2\sqrt{\lambda}\, f_1 &= i\lambda (\U1-i \U3) + i (1+z_{12}^3-i z_{21}^2)
(\U5-i\U6)\label{reveq2}.
\end{align} 

However, it follows from \eqref{matb} that $-i= (f_1,g_1)$, so orthonormality of $\cal U$ requires that
\begin{align*}
-4 \theta_2 \lambda i &= 2 i\lambda  (z_{12}^3-1 -i z_{21}^2)- 2 i\lambda (1+z_{12}^3 -i
z_{21}^2)\\
&=-4 i \lambda,
\end{align*}
so that $\theta_2=1$.
Therefore, in order to construct $\cal U$, it is necessary to
determine real-valued functions $z_{12}^3(t,z)$, $z_{21}^2(t,z)$, $\lambda(t,z)$, and orthonormal vector fields
$\U1(t,z)$, $\U3(t,z)$, $\U5(t,z)$, $\U6(t,z)$ in $\R^6$ satisfying \eqref{reveq1} and 
\eqref{reveq2} (with $\theta_2=1$ in this latter equation). 
However, as $(f_1,\bar f_1) =e^{\omega}$ and $(g_1,\bar g_1) =
e^{\omega^+}$ we see that
\begin{align*}
&2\lambda e^{\omega^+} = (z_{12}^3-1)^2 +(z_{21}^2)^2 +  \lambda^2,\\
&2\lambda e^{\omega} = (1+z_{12}^3)^2 +(z_{21}^2)^2 +  \lambda^2.\end{align*}
Thus
\begin{equation}
z_{12}^3 = \tfrac{1}{2} \lambda (e^{\omega} -e^{\omega^+}),\label{z123}
\end{equation}
and
\begin{equation*}
(z_{21}^2)^2 + \left(1- \tfrac{1}{2} \lambda (e^{\omega} +e^{\omega^+})\right)^2= \lambda^2 (e^{\omega
+\omega^+}-1).
\end{equation*}

As $\lambda$ and $(e^{\omega
+\omega^+}-1)$ are both positive, we may define $z_{21}^2$ and $\lambda$ by taking
\begin{align}
&\lambda =\tfrac{2}{e^{\omega} + e^{\omega^+}+2\cos t  \sqrt{ e^{\omega
+\omega^+}-1}},\label{lambda}\\
&z_{21}^2 = \lambda \sin t \sqrt{ e^{\omega
+\omega^+}-1}. \label{z212}
\end{align}

We restrict  $t$ to lie on  a suitable subinterval $I$ of $(0,\pi)$, in order to ensure that $\lambda$ is well
defined and $z_{21}^2>0$.  

We have now obtained, through \eqref{z123},  \eqref{lambda} and \eqref{z212}, formulae for $z_{12}^3$, $\lambda$ and $z_{21}^2$ in terms of 
 $\omega$, $\omega^+$ and $\gamma^+$.  We next obtain $\U1$,  $\U3$, $\U5$ and $\U6$ as the solutions to \eqref{reveq1} and \eqref{reveq2}. In the next section we will discuss a particular special case of the construction detailed in this section, so we will write down the  solution to \eqref{reveq1} and \eqref{reveq2} and verify their properties explicitly (for which we used Mathematica), although the properties we obtain may be deduced directly from  \eqref{reveq1} and \eqref{reveq2}.

So, if 
\begin{equation*}
C= \sqrt{\tfrac{\lambda}{e^{\omega
+\omega^+}-1}}\ ,
\end{equation*}
we find from \eqref{reveq1} and \eqref{reveq2} that
\begin{align}
&\U1+i \U3=   -C\left((\sqrt{ e^{\omega
+\omega^+}-1} +e^{-it+\omega})g_1 + i e^{-it} \bar f_1\right), \label{U13}\\
&\U5+i \U6= C \left(e^{-it}  g_1 +i  (\sqrt{ e^{\omega
+\omega^+}-1}+ e^{-it+\omega^+}) \bar f_1\right).\label{U56}
\end{align} 
It is now a straightforward computation using \eqref{matb} to verify that $\U1$ up to $\U6$ defined by 
\eqref{U24},\eqref{U13} and  \eqref{U56}
are orthonormal vectors and that, using \eqref{volume2},
\begin{align*}
\hbox{vol}(\U1,\dots,\U6) &= -\tfrac{1}{e^{\omega+\omega^+}-1} \hbox{vol}(f_0,f_1,\bar f_1,\bar g_1,g_1,g_0)\\
&=1.
\end{align*}
Thus ${\cal U}=(\U1,\ldots ,\U6):M\to SO(6)$, so that
$\Omega={\mathcal U}^{-1}d{\mathcal U}$ is a skew symmetric matrix whose second and fourth columns have the correct form.

It remains to find the linearly independent 1-forms $\omega_1,\omega_2,\omega_3$  on $M$, and real valued  functions $z_{22}^3$, $z_{32}^3$, $a$ and $b$ on $M$ such that the entries of $\Omega=\mathcal U^{-1}d\mathcal U$ are as 
given in (33) of \cite{bovr}. 

Since $d\U2(\partial/\partial t)=0$, for (33) of \cite{bovr} to hold we need that $\omega_2(\partial/\partial t)=\omega_3(\partial/\partial t)=0$. Also, from \eqref{XiY} with $\theta_1=1$, we have
\begin{align*}
d\U2(\e2)&=\sqrt\lambda\, dg(\partial/\partial x) = \sqrt\lambda\, d\U2(\partial/\partial x),\\
d\U2(\e3)&=\sqrt\lambda\, dg(\partial/\partial y) = \sqrt\lambda \,d\U2(\partial/\partial y),
\end{align*}
from which it follows that $\e2-\sqrt\lambda \,\partial/\partial x$ and $\e3-\sqrt\lambda \,\partial/\partial y$ are multiples of $\partial/\partial t$.  
Thus 
\begin{equation}
\omega_2=\tfrac{1}{\sqrt\lambda}dx, \quad \omega_3=\tfrac{1}{\sqrt\lambda}dy,\label{o23}
\end{equation}
which also ensures that \eqref{reveq2} holds with $\theta_2=1$. 

We now consider the columns of $\Omega$ other than the second and fourth.  These have the correct form if and only if we have that, 
modulo $\U2$ and $\U4$,
\begin{align} 
\begin{split}
d(\U1+i \U3) &\equiv i \left((1+z_{12}^3)\w1 + z_{22}^3 \w2 +  z_{32}^3\w3\right) (\U1+i \U3)\\ 
&\qquad + (c \lambda^{-\tfrac{1}{2}} dz - i\lambda \w1) (\U5+i \U6), \quad (\hbox{mod }\U2,\U4), \label{dU13}
\end{split}\\
\begin{split}
d(\U5+i \U6) &\equiv -(\bar c \lambda^{-\tfrac{1}{2}} d\bar z + i\lambda \w1)(\U1+i\U3)\\
&\qquad -i \left((z_{12}^3-1)\w1 +  z_{22}^3 \w2+ z_{32}^3 \w3\right)(\U5+i\U6), \quad (\hbox{mod }\U2,\U4),\label{dU56}
\end{split}
\end{align}
where $c = -b-i a$. 
In particular, $\w1$ must satisfy
\begin{equation}\label{o1}
 \big(d(\U1+i \U3),\U1-i\U3\big) +\big(d(\U5+i \U6),\U5-i\U6\big)= 4 i \w1.
 \end{equation}
Computing the lefthand side explicitly using \eqref{U13} and \eqref{U56}, we see that \eqref{o1} holds if and only if 
\begin{align}
&\w1(\tfrac{\partial}{\partial t})= -\tfrac{1}{2},\\
&\w1(\bar \partial)=-\tfrac{1}{4} i \tfrac{2 i \,\overline {\gamma^+} + e^{\omega
+\omega^+} \bar \partial (\omega-\omega^+)}{e^{\omega
+\omega^+}-1}.\label{o1dbar}
\end{align}
We use these expressions to define $\omega_1$, in which case \eqref{o23}, implies that $\w1$, $\w2$ and $\w3$ are linearly independent $1$-forms on $M$. 

A straightforward computation using \eqref{U13} and \eqref{U56} now shows that  
\begin{align*}
 &\big(d(\U1+i \U3),\U5-i\U6\big)(\tfrac{\partial}{\partial t})= i \lambda,\\
 &\big(d(\U1+i \U3),\U1-i\U3\big)(\tfrac{\partial}{\partial t})-\big(d(\U5+i \U6),\U5-i\U6\big)
 (\tfrac{\partial}{\partial t}) + 2 i z_{12}^3 =0,\\
 &\big(d(\U1+i \U3), \U5-i\U6\big)(\bar \partial) =-\tfrac{1}{2} \lambda \tfrac{2 i  \,\overline {\gamma^+} + e^{\omega
+\omega^+} \bar \partial (\omega-\omega^+)}{e^{\omega
+\omega^+}-1}=-2i\lambda\omega_1(\bar \partial),
 \end{align*}
so it only remains to define the complex-valued function $c$ on $M$ and the real
valued functions $z_{22}^3$ and $z_{32}^3$ in such a way that \eqref{dU13} and  \eqref{dU56} hold. This may be done explicitly and uniquely 
by calculating  $\big(d(\U1+i \U3), \U5-i\U6\big)(\partial)$,
$\big(d(\U1+i \U3), \U1-i\U3\big)(\tfrac{\partial}{\partial x})$, and $\big(d(\U1+i \U3),
\U1-i\U3\big)(\tfrac{\partial}{\partial y})$. We have thus proved the following theorem.
\begin{theorem}  Let $f:S \rightarrow S^5(1)$ be a minimal immersion with non-degenerate non-circular ellipse of curvature at every point. Then there exists a minimal Lagrangian submanifold of $CP^3(4)$, admitting a foliation by asymptotic curves, for which the construction described in \cite{bovr} yields $f$ and its $(+)$transform. 
\end{theorem}
  
We remark that, since we have shown in the previous section that such minimal surfaces are part of a sequence, minimal Lagrangian submanifolds of type 3 also form a sequence. However, up to now, we
do not know geometrically (without using this detour over  minimal surfaces) how to associate one with its successor.   

In the next section we will give an example in which we can explicitly describe 
the reverse construction detailed in this section.

\section{Lawson's bipolar surfaces}
Let $f:S\to S^5(1)$ be a minimal immersion with non-degenerate non-circular ellipse of curvature, and let $z$ be an adapted complex coordinate for $f$. Using the notation of the previous section, we will consider the special case in which the invariant $\gamma^+=(\partial f_1,g_1)$ is identically zero. We will show that, in this case, $f$ is the bipolar surface in the sense of Lawson \cite{la} of a minimal surface in $S^3(1)$.

We begin by noting that if $\gamma^+=0$ then \eqref{eq24} implies that $\omega= \omega^+$ while \eqref{eq24} and Lemma \ref{lem1} impliy that 
$\omega$ is a positive solution  of the following partial
differential equation:
\begin{equation}
\label{sinh1}
 \partial\bar{\partial}\omega=-2 \sinh \omega
 +\tfrac{1}{e^{2\omega}-1} \vert \partial
 \omega\vert^2.
\end{equation}   
Conversely given a positive solution of the above differential equation, there exists a corresponding
minimal surface in $S^5(1)$ with non-circular non-degenerate ellipse of curvature and induced metric $2e^\omega|dz|^2$. 

It is convenient to rewrite the above differential equation by making the substitution $e^{\omega} = \cosh \eta$. A short calculation shows that $\omega$ satisfies \eqref{sinh1} if and only if 
the function $\eta$ satisfies the sinh-Gordon equation
\begin{equation}
\label{sinh2}
 \partial\bar{\partial}\eta=- \sinh \eta.
\end{equation}
We recall \cite{fpps} that a solution $\eta$ of the sinh-Gordon equation determines an $S^1$-family of  non-totally geodesic  minimal immersions in $S^3(1)$ whose induced metric is $e^\eta|dz|^2$, and we will see that $f$ is the bipolar in $S^5(1)$ of the minimal immersion in this family for which the coordinate curves are the lines of curvature.

Specialising the formulae of the previous section to the case $\gamma^+=0$, we obtain
\begin{equation}\label{invF} \begin{aligned}
&\lambda =\tfrac{1}{\cosh \eta + \cos t \sinh \eta},\\
&z_{21}^2= \tfrac{\sin t \sinh \eta}{\cosh \eta + \cos t \sinh \eta},\\
&z_{12}^3 = 0,\\
&\w1=-\tfrac{1}{2} dt,\\
&\w2=\tfrac{1}{\sqrt{\lambda}} dx,\\
&\w3=\tfrac{1}{\sqrt{\lambda}} dy,\\
&b=\tfrac{1}{2} \lambda^{3/2} \eta_y  \sin t,\\
&a=\tfrac{1}{2} \lambda^{3/2} \eta_x  \sin t,\\
&z_{32}^3=\tfrac{1}{2} \lambda^{3/2} \eta_x (\cos t \cosh \eta +\sinh \eta),\\
&z_{22}^3=-\tfrac{1}{2} \lambda^{3/2} \eta_y (\cos t \cosh \eta +\sinh \eta).
\end{aligned}
\end{equation}
In particular,
\begin{align*}
&\partial/\partial t= -\tfrac{1}{2} \e1,\\
&\partial/\partial x=\tfrac{1}{\sqrt{\lambda}} \e2,\\
&\partial/\partial y =\tfrac{1}{\sqrt{\lambda}} \e3.
\end{align*} 

Substituting these expressions into (25) of \cite{bovr} we see that the horizontal lift $F$
to $S^7(1)$ of the minimal Lagrangian immersion into $\C P^3(4)$ corresponding to $f$ satisfies the following system of differential equations:
\begin{align}
&F_{tt}=-F/4,\label{Ftt}\\
&F_{tx}=\tfrac{-\left( \imag  + \sin t\,\sinh \eta \right) }
   {2\,\left( \cosh \eta + \cos t\,\sinh \eta \right) }F_x,\label{Ftx}\\
&F_{ty}=\tfrac{\imag  - \sin t\,\sinh \eta}
   {2(\cosh \eta + \cos t\,\sinh \eta)}F_y,\label{Fty}\\
\begin{split}&F_{xx}=-(\cosh \eta +\cos t\,\sinh \eta)F +2(\sin t\,\sinh \eta\,-i)F_t\\& \qquad\qquad
+ \tfrac{\eta_x\,\left( \cos t\,\cosh \eta + \imag \,\sin t + 
       \sinh \eta \right) }{2\,\left( \cosh \eta + \cos t\,\sinh \eta \right) }F_x-\tfrac{\eta_y\,\left( \cos t\,\cosh \eta - 
         \imag \,\sin t + \sinh \eta \right)}{2\,
     \left( \cosh \eta + \cos t\,\sinh \eta \right) }F_y,
\end{split}\label{Fxx}\\
&F_{xy}=\tfrac{\eta_y \,\left( \cos t\,\cosh \eta + \imag \,\sin t + 
       \sinh \eta \right) }{2\,\left( \cosh \eta + \cos t\,\sinh \eta \right) }F_x+
  \tfrac{\eta_x \,\left( \cos t\,\cosh \eta - \imag \,\sin t + 
       \sinh \eta \right) }{2\,\left( \cosh \eta + \cos t\,\sinh \eta \right) } F_y,\label{Fxy}\\
\begin{split}&F_{yy}=-(\cosh \eta +\cos t\,\sinh \eta)F+2(\imag  + \sin t\,\sinh \eta)F_t\\& 
\qquad\qquad -
  \tfrac{\eta_x\,\left( \cos t\,\cosh \eta + 
         \imag \,\sin t + \sinh \eta  \right) }{2\,
     \left( \cosh \eta + \cos t\,\sinh \eta \right) }F_x+
  \tfrac{\eta_y\,\left( \cos t\,\cosh \eta - \imag \,\sin t + 
       \sinh \eta \right) }{2\,\left( \cosh \eta + \cos t\,\sinh \eta \right) }F_y.\end{split}\label{Fyy}
\end{align}
It follows from \eqref{Ftt} that we may write
\begin{equation}
F(t,u,v) = G_1(x,y) \cos \tfrac{t}{2} + i G_2(x,y) \sin \tfrac{t}{2},\label{FGG}
\end{equation}
for suitable  $\mathbb C^4$-valued functions $G_1$ and $G_2$. 
Substituting this into \eqref{Ftx} and carrying out significant but elementary simplification
we find that
\begin{equation}
(G_2)_x = -e^{-\eta} (G_1)_x, \label{Gx}
\end{equation}
while similar reasoning using \eqref{Fty} gives that 
\begin{equation}
(G_2)_y = e^{-\eta} (G_1)_y. \label{Gy}
\end{equation}
Using \eqref{Gx} and \eqref{Gy} we find,  after some calculation,  that \eqref{Fxx}, \eqref{Fxy} and \eqref{Fyy} are equivalent to 
\begin{align}
\begin{split}
(G_1)_{xx}& =\tfrac{1}{2} \eta_x (G_1)_x - \tfrac{1}{2} \eta_y (G_1)_y +G_2-e^\eta G_1,\\
(G_1)_{xy}& =\tfrac{1}{2} \eta_y (G_1)_x + \tfrac{1}{2} \eta_x (G_1)_y,\label{G1min}\\ 
(G_1)_{yy}& =-\tfrac{1}{2} \eta_x (G_1)_x + \tfrac{1}{2} \eta_y (G_1)_y -G_2-e^\eta G_1.
\end{split}
\end{align}

We now note that since $|F|=1$, \eqref{FGG} implies that $|G_1|=|G_2|=1$ and that $G_1$ is real  orthogonal to $G_2$.  The horizontality of $F$ further shows that $G_1$ and $G_2$ are unitarily orthogonal.

The coefficients in the system \eqref{G1min} are all real, so the real subspace spanned by $G_1$, $G_2$, $(G_1)_x$ and
$(G_1)_y$ is constant. We identify this subspace with $\mathbb R^4$ by picking 
an orthonormal basis.
But the system \eqref{G1min} is exactly that of a minimal surface $G_1$ in $S^3(1)$ with unit normal  $G_2$, induced metric
$ds^2 =e^\eta \vert dz\vert^2$ with the complex coordinate chosen such that the second fundamental form $\tilde{II}$ of  $G_1$ satisfies 
$(\tilde{II}(\partial,\partial), \tilde{II}(\partial,\partial))=1/4$. In particular, since $\tilde{II}(\partial/\partial x,\partial/\partial y)=0$, we see that the coordinate curves of $G_1$ are the lines of curvature.

Applying now the definition of $U_4$ of \cite{bovr} and the expressions for $\omega_1$, $\omega_2$ and $\omega_3$ obtained in \eqref{invF}, we get that 
\begin{align*}
f& = 
\tfrac{1}{i\sqrt{2}}(F \wedge (-2 F_t)-\lambda F_x \wedge F_y)\subset \Lambda^2{\mathbb C}^4.
\end{align*}
An easy calculation using \eqref{FGG} now shows that 
\begin{equation}\label{bip}
f=\tfrac{1}{\sqrt{2}} (ie^{-\eta} G_{1x} \wedge G_{1y}-G_1\wedge G_2).
\end{equation}

According to Lawson in \cite{la}, the {\em bipolar surface} of the minimal surface $G_1$ in $S^3(1)$ with unit normal $G_2$ is the surface in $S^5(1)$ given by $G_1\wedge G_2$ in $\Lambda^2 \R^4 = \R^6$. This is a minimal surface in $S^5(1)$. If we include $\Lambda^2 \R^4$ in  $\Lambda^2 \C^4$ via $v\mapsto (1/\sqrt 2)(v-i\star v)$, where $\star$ denotes the Hodge star operator on  $\Lambda^2 \R^4$, then it follows immediately from \eqref{bip} that $f$ is the bipolar of $G_1$. 

Conversely, let $G_1(z)$ be a non-totally geodesic minimal immersion in $S^3(1)$ with unit normal $G_2(z)$. We may assume that $z$ is such that the induced metric is $ds^2 =e^\eta \vert dz\vert^2$, where $\eta$ satisfies \eqref{sinh2}  and that the second fundamental form  $\tilde{II}$ of  $G_1$ satisfies 
$(\tilde{II}(\partial,\partial), \tilde{II}(\partial,\partial))=1/4$. Then $G_1$ and $G_2$ will satisfy the system \eqref{G1min}, so if we define $F$ using \eqref{FGG} then $F$ is  horizontal in $S^7(1)$ and satisfies the system \eqref{Ftt}-\eqref{Fyy}.  In particular, we have that $z_{12}^3=0$ and $\omega_1(\bar\partial)=0$. If we apply the construction of \cite{bovr} to the projection of $F$ to $\C P^3$, then we will obtain the bipolar $f(z)$ of $G_1(z)$. It then follows from \eqref{z123} that $\omega=\omega^+$, so that \eqref{o1dbar} gives that $\gamma^+=0$.

We have therefore proved the following theorem. 

\begin{theorem} Let $f: S \rightarrow S^5(1)$ be a minimal immersion with non-degenerate 
non-circular ellipse of curvature. Then the following three statements are equivalent:
\begin{enumerate} 
\item $\gamma^+=0$, 
\item $f$ is the bipolar surface of a non-totally geodesic minimal surface in $S^3(1)$,
\item the $(+)$transform $f^+$ is the reflection of $f$ in a great subsphere of $S^5(1)$.
\end{enumerate}
\end{theorem}

\noindent {\em J. Bolton, Dept of Mathematical Sciences, University of Durham,

Durham DH1 3LE, UK. \quad  E-mail: john.bolton@dur.ac.uk

\noindent L. Vrancken,   LAMATH, ISTV2, Universit\'e de Valenciennes, 
Campus du Mont Houy, 59313 Valenciennes Cedex 9, France. \quad E-mail: 
luc.vrancken@univ-valenciennes.fr}  


\begin{thebibliography}{1}

\bibitem{boscvrwo}
J.~Bolton, C.~Scharlach, L.~Vrancken and L.~M.~Woodward.
\newblock From certain minimal Lagrangian submanifolds of the 
3-dimensional complex 
projective space to minimal surfaces in the 5-sphere.
\newblock {\em Proceedings of the Fifth Pacific Rim
Geometry Conference, Tohoku University}. Tohoku Mathematical Publication No. {\bf 20} 23--31, 2001.

\bibitem{boscvr}
J.~Bolton,  C.~Scharlach and L.~Vrancken. 
\newblock From  surfaces in the 5-sphere to 3-manifolds in complex 
projective 3-space.
\newblock {\em Bull. Austral. Math. Soc.} {\bf 66} 465--475, 2002.


\bibitem{bovr}
J.~Bolton and  L.~Vrancken.
\newblock Ruled minimal Lagrangian submanifolds of complex projective 3-space.
\newblock To appear in {\em Asian J. Math.}

\bibitem{bovrwo}
J.~Bolton, L.~Vrancken and L.M.~Woodward.
\newblock Totally real minimal surfaces with non-circular ellipse of curvature in the nearly K\"ahler $S^6$.
\newblock {\em J. London Math. Soc (2)} {\bf 56} 625--644, 1997. 

\bibitem{bowo}
J.~Bolton and  L.M.~Woodward.
\newblock Congruence theorems for harmonic maps from a Riemann surface into $\C P^n$ and $S^n$.
\newblock {\em J. London Math. Soc (2)} {\bf 45} 363--376, 1992. 

\bibitem{br}
R.L.~Bryant.
\newblock Second order families of special Lagrangian 3-folds.
\newblock Preprint.

\bibitem{cavr}
I. Castro and L. Vrancken. 
\newblock Minimal Lagrangian 
submanifolds in $\mathbb C {P}\sp 3$ and the sinh-Gordon equation. Dedicated 
to Shiing-Shen Chern on his 90th birthday.  
\newblock {\em Results Math.}  {\bf 40 no. 
1-4} 130--143, 2001. 

\bibitem{ch}
B.-Y.~Chen, F.~Dillen, L.~Verstraelen and L.~Vrancken.
\newblock Totally real submanifolds of ${\mathbb C}P^n$ satisfying a basic
equality.  
\newblock {\em Arch. Math.} {\bf 63} 553--564, 1994.

\bibitem{fpps} D.~Ferus, F.~Pedit, U.~Pinkall and I.~Sterling.
\newblock Minimal tori in $S^4$.
\newblock {\em J. reine angew. Math.} {\bf 429} 1--47, 1992.

\bibitem{la}
H.B. Lawson. 
\newblock Complete minimal surfaces in $S^3$. 
\newblock {\em Annals of Mathematics} {\bf 92} (1970) 335-374.

\end{thebibliography}
\end{document}